\documentclass{CSML}

\def\dOi{13(1:5)2017}
\lmcsheading%
{\dOi}
{1--6}
{}
{}
{Apr.~28, 2016}
{Feb.~\phantom02, 2017}
{}

\ACMCCS{[{\bf Theory of computation}]: ???}

\usepackage{hyperref}

\newcommand\nn{{\{0,1\}^n}}
\newcommand\gn{{\{0,1\}^{{n \choose 2}}}}
\newcommand\at{ \{0, 1\}^t }
\newcommand\as{ \{0, 1\}^s }
\newcommand\rr{{\bf r}}
\newcommand\bs{{\bf s}}
\newcommand\bp{{\bf p}}
\newcommand\bq{{\bf q}}
\newcommand\pro{(G, {\sf lab}, F, S)}

\begin{document}

\title[A feasible interpolation for random resolution]
      {A feasible interpolation for random resolution}

\author[J.~Kraj\'{\i}\v{c}ek]{Jan Kraj\'{\i}\v{c}ek}
\address{Faculty of Mathematics and Physics,
Charles University in Prague}
\email{krajicek@karlin.mff.cuni.cz}

\begin{abstract}
Random resolution, defined by Buss, Kolodziejczyk and Thapen (JSL, 2014), is a sound propositional
proof system that extends the resolution proof system by the possibility to augment any set of initial
clauses by a set of randomly chosen clauses (modulo a technical condition).
We show how to apply the general feasible interpolation theorem for semantic derivations
of Kraj\'{\i}\v cek (JSL, 1997) to random resolution. As a consequence
we get a lower bound for random resolution refutations of the clique-coloring
formulas.
\end{abstract}

\maketitle

\section{Introduction}
Assume $A_1, \dots, A_m, B_1, \dots, B_\ell$ is an unsatisfiable set of clauses
in variables partitioned into three disjoint sets 
$\bp$, $\bq$ and $\rr$, with clauses $A_i$ containing only variables 
from $\bp$ and $\bq$ while clauses
$B_j$ contain only variables from $\bp$ or $\rr$.

Feasible interpolation for resolution \cite[Thm.6.1]{Kra-interpol} says 
that if the set has a resolution refutation
with $k$ clauses then there is a circuit
of size\footnote{The bound $k n^{O(1)}$ is derived from a general interpolation theorem for semantic derivations
whose framework we also use below; a bit better bound (proportional to the size of the refutation
and hence $O(k n)$) can be proved by resolution specific arguments.} $k n^{O(1)}$, where $n$ is the number of variables $\bp$,
with inputs $\bp$ that outputs $1$ on all $\bp := {\bf a} \in \nn$
for which $\bigwedge_i A_i({\bf a}, \bq)$ is satisfiable and $0$ on all $\bf a$
for which
$\bigwedge_j B_j({\bf a}, \rr)$ is satisfiable. Moreover, if variables $\bp$ occur only
positively in clauses $A_i$ 
then the interpolating circuit can be required to be monotone.

The monotone version
can then be applied to the clique-coloring clauses \cite[Def.7.1]{Kra-interpol}
where there are ${n \choose 2}$ variables $\bp$ indexed by unordered pairs $i, j$ of different elements
from $[n]:=\{1, \dots, n\}$, $\omega \cdot n$ variables $\bq$ indexed by elements of
$[\omega] \times [n]$ and $n \cdot \xi$ variables $\rr$ indexed by elements of $[n]\times [\xi]$,
with $n \geq \omega > \xi \geq 1$: 

\begin{enumerate}

\item $\{q_{u 1}, \dots, q_{u n}\}$, for each $u \in [\omega]$

\item $\{\neg q_{u i}, \neg q_{v i}\}$, for $u \neq v \in [\omega]$ and $i \in [n]$

\item $\{\neg q_{u i}, \neg q_{v j}, p_{i j}\}$, for $u \neq v \in [\omega]$ and $i \neq j \in [n]$

\item $\{r_{i 1}, \dots, r_{i \xi}\}$, for each $i \in [n]$

\item $\{\neg r_{i u}, \neg r_{i v}\}$, for each $u \neq v  \in [\xi]$ and $i \in [n]$

\item $\{\neg r_{i v}, \neg r_{j v}, \neg p_{i j}\}$, for $v \in [\xi]$ and $i \neq j \in [n]$

\end{enumerate}
The clauses in the first three items comprise the set $Clique_{n,\omega}$ and the clauses in the
last three items comprise the set $Color_{n,\xi}$, They have only variables $\bp$ in
common and these occur only positively in $Clique_{n,\omega}$. The assignments ${\bf a}$ to
$\bp$ for which $Clique{_n,\omega}({\bf a}, \bq)$ is satisfiable can be identified with 
undirected graphs on $[n]$ without loops and 
having a clique of size at least $\omega$ while those $\bf a$
for which $Color_{n,\xi}({\bf a}, \rr)$ is satisfiable are $\xi$-colorable graphs.
Hence $Clique_{n,\omega} \cup Color_{n, \xi}$ is unsatisfiable as $\xi < \omega$ and the monotone feasible
interpolation combined with the Alon-Boppana \cite{AloBop} exponential lower for monotone circuits
separating the two classes of graphs implies that all resolution refutations of the set 
must have an exponential number of clauses, cf.\cite[Sec.7]{Kra-interpol}.

Buss, Kolodziejczyk and Thapen \cite[Sec.5.2]{BKT} defined the notion
of $\delta$-random resolution (the definition is attributed in \cite{BKT}
to S.~Dantchev). The motivation for introducing the proof system
came from bounded arithmetic; the proof system simulates an interesting theory.
A {\em $\delta$-random resolution refutation distribution} of
a set of clauses $\Psi$ (\cite{BKT} considers only narrow clauses because of the
specific problem studied there) is a random distribution $(\pi_{\bs}, \Delta_{\bs})_{\bs}$
such that $\pi_\bs$ is a
resolution refutation of $\Psi \cup \Delta_{\bs}$, and where the following technical condition is
satisfied:
\begin{itemize}
\item  
any fixed truth assignment to all variables
satisfies the set of clauses $\Delta_{\bs}$
with probability at least $1 - \delta$. 
\end{itemize}
The number of clauses in such a random refutation is the maximal number of
clauses among all $\pi_{\bs}$. 
Note that it is a sound proof system in the sense that any refutable set $\Psi$
is indeed unsatisfiable: if ${\bf a}$ would be a satisfying assignment for $\Psi$ then,
by the condition above, ${\bf a}$ would satisfy also some $\Delta_\bs$ and hence $\pi_\bs$
would be a resolution refutation of a satisfiable set of clauses which is impossible.
Variants of the definition of this proof system and its properties 
are studied in \cite{PudTha}. 

The presence of the clauses $\Delta_\bs$ spoils the separation of the $\bq$ and $\rr$
variables in initial clauses
and this seems to prohibit any application of the feasible interpolation 
method. The point of this note is to show that, in fact, the construction behind the
general feasible interpolation 
theorem  \cite{Kra-interpol} for semantic derivations based on communication complexity
does apply here fairly straightforwardly.

We recall some feasible interpolation preliminaries from \cite{Kra-interpol} in
Section \ref{1}. In Section \ref{2} we prove monotone feasible interpolation for random resolution 
and this will yield the following lower bound for
random resolution refutations of the clique-coloring clauses.

\begin{thm}\label{main}
Let $n \geq \omega > \xi \geq 1$ and $\xi^{1/2}\omega \le 8n/\log n$. 
Assume $\delta < 1$ and let $(\pi_{\bs}, \Delta_{\bs})_{\bs}$
be a $\delta$-random resolution refutation distribution of 
$Clique_{n,\omega} \cup Color_{n,\xi}$ with $k$ clauses.
Put $d := \max_\bs |\Delta_\bs|$.

Then:

\begin{enumerate}

\item If $d \delta < 1$ then $k \geq (1 - d\delta^{1/2})n^{\Omega(\xi^{1/2})}$. 

\item $k \geq \min(1/(2 \delta^{1/2}), n^{\Omega(\xi^{1/2})})$. 

\end{enumerate}

\end{thm}
The proof of this theorem will be given at the end of Section \ref{2}.
We only remark that for tree-like refutations a feasible interpolation
via ordinary randomized Karchmer-Wigderson protocols 
follows from \cite{Kra-interpol}
immediately and it yields an exponential lower bound for formulas formalizing Hall's
theorem as described in \cite[Sec.4]{Kra-game}.

We will give below a detailed formulation of constructions from 
\cite{Kra-interpol} needed here but we will
not repeat the arguments 
from that paper. For more general background on proof complexity 
the reader may consult \cite{kniha,Pud-survey}.

\section{Feasible interpolation via protocols} \label{1}

We review the needed material from \cite{Kra-interpol} just for the case of monotone
interpolation and the clique-coloring clauses (but it is quite representative). 
Identify undirected graphs without loops on $[n]$ with strings from $\gn$.
Note that indices of $\bp$ variables correspond to pairs of different
vertices and hence the truth value an assignment $\bf a$ gives to a particular
$\bp$-variable indicates whether or not 
the edge corresponding to the variable is in the graph $\bf a$.

Let $U \subseteq \gn$ be the set of graphs having a clique of size at least $\omega$
and let $V \subseteq \gn$ be the set of $\xi$-colorable graphs. Let the monotone Karchmer-
Wigderson function $KW^m(u,v)$
be a multi-function defined on $U \times V$ whose valid value on
a pair $(u,v) \in U \times V$ is any edge 
(i.e. unordered pair $i\neq j \in [n]$) that is present in $u$ but not in $v$.

The method in \cite{Kra-interpol} 
extracts from a resolution refutation of $Clique_{n,\omega} \cup Color_{n,\xi}$
a protocol for a communication between two players, one holding $u$ and the other one $v$,
who want to find a valid value for $KW^m(u,v)$. The protocols in \cite{Kra-interpol} are,
however, more complex than just binary trees as in the ordinary communication
complexity set-up of \cite{KW}.
 
A monotone protocol for computing $KW^m$ in the sense of  \cite[Def.2.2]{Kra-interpol}
is a 4-tuple $\pro$
satisfying the following conditions:
\begin{enumerate}

\item $G$ is a directed acyclic graph that has one root (the in-degree $0$ node)
denoted $\emptyset$.

\item The nodes with the out-degree $0$ are leaves
and they are labelled by the mapping ${\sf lab}$. The mapping ${\sf lab}$
assigns an
element of $[{n \choose 2}]$ (i.e., a potential edge) to each leaf in $G$.

\item $S(u,v,x)$ is a function (called the strategy)
that assigns to a node $x \in G$ 
and a pair $u \in U$ and $v \in V$ a 
node $S(u,v,x)$ reachable from the node $x$ by one edge.

\item For every 
$u \in U $ and $v \in V$, $F(u,v) \subseteq
G$ is a set (called the consistency condition) satisfying:
\begin{enumerate}

\item $\emptyset \in F(u,v)$,

\item $x \in F(u,v) \longrightarrow S(u,v,x) \in F(u,v)$,

\item if $x \in F(u,v)$ is a leaf and ${\sf lab}(x)= \{i,j\}$, then
$u_{i,j} = 1 \wedge v_{i,j} = 0$ holds.
\end{enumerate}
\end{enumerate}\medskip

\noindent The {\em size} of $\pro$ is the cardinality of $G$ 
and its {\em communication complexity} is the minimal $t$ such that  for every $x \in G$
the communication complexity for 
the players (one knowing $u$ and $x$, the other one $v$ and $x$)
to decide $x \in_? F(u,v)$ or to compute $S(u,v,x)$
is at most $t$.

Put $s := n \cdot \omega$ and identify strings from $\as$ with assignments to 
$\bq$-variables, 
and similarly put $t := n \cdot \xi$ and identify strings from $\at$ with assignments to
$\rr$-variables. For any $u \in U$ fix $q^u \in \as$ such that $(u, q^u)$ satisfies
all clauses from $Clique_{n,\omega}$ and for $v \in V$ fix $r^v \in \at$ such that
$(v, r^v)$ satisfies all clauses of $Color_{n,\xi}$.

The protocol $\pro$ for $KW^m$
constructed in \cite[Thm.5.1 and Thm.6.1]{Kra-interpol}
from a resolution refutation $\pi$ of $Clique_{n,\omega} \cup Color_{n,\xi}$
having $k$ steps has $k + {{n}\choose{2}}$ nodes:
$k$ nodes corresponding to the clauses of $\pi$ are the inner nodes
and ${n \choose 2}$ other nodes
are the leaves and these are labelled by the ${n \choose 2}$ possible values of the 
multi-function $KW^m$. The consistency condition $x \in F(u,v)$ for a node $x$ corresponding to 
a clause $C$ of $\pi$ is defined by the condition
that the assignment $(v, q^u, r^v)$ falsifies $C$, and for a leaf by the 
condition that the label is a valid value of $KW^m$
for the pair $(u,v)$. The strategy $S$ (whose exact definition
we do not need) navigates from the root (the end-clause of $\pi$) through $\pi$
towards the initial clauses and the construction shows that sooner or later it encounters
a situation that allows it to compute a valid value of $KW^m$ and move to the leaf with the appropriate
label. The construction is fairly general and we shall formulate in Theorem \ref{3.2}
its one particular feature.

For a set $\Delta$ of clauses in variables $\bp, \bq$ and $\rr$ define a multifunction $F_\Delta$
on $U \times V$ whose valid value on a pair $(u,v)$ is any valid value of $KW^m(u,v)$
and also a new value $\bot$ {\em provided} that $(v, q^u, r^v)$ falsifies some clause in $\Delta$.
Note the similarity of the condition permitting the value $\bot$ with the consistency
condition in the protocol just discussed. 

Now we recall a particular fact about the existence of
protocols provided by the constructions in the proofs of \cite[Thm.5.1 and Thm.6.1]{Kra-interpol}
(again we restrict ourselves to the clique-coloring formulas and the monotone case).

\begin{thm}[{\cite{Kra-interpol}}] \label{3.2}
Assume that $\Delta$ is a set
of clauses in variables $\bp, \bq$ and $\rr$ and that
$\pi$ is a resolution refutation of the set $Clique_{n,\omega} \cup Color_{n,\xi} \cup \Delta$
and that $\pi$ has $k$ steps.

Then there is a protocol $\pro$ 
for $F_\Delta$ of size $k + {n \choose 2}$ whose strategy has the communication complexity 
at most $2 + 2\log n$ and whose consistency condition has the communication complexity $2$.

Further, the existence of a 
protocol for $KW^m$ on $U'\times V' \subseteq U \times V$ of size
$k'$ and monotone communication complexity $O(\log n)$ 
implies the existence of a monotone circuit
of size at most $k' \cdot n^{O(1)}$ separating $U'$ from $V'$.

\end{thm}

The part about the existence of a circuit is in \cite{Kra-interpol} proved using
a result from \cite{Raz95}; a stand alone proof can be found in 
\cite[Sec.2.4]{Kra-ds1}.

\section{The lower bound} \label{2}

For $(u,v) \in U \times V$ define $w(u,v) := (v, q^u, r^v)$ and 
for $X \subseteq U$ and $Y \subseteq V$ define $W(X,Y) \subseteq \gn \times \as \times \at$ 
to be the set of all tuples $w(u,v)$ for $(u,v) \in X \times Y$.

Assume $(\pi_{\bs}, \Delta_{\bs})_{\bs}$ is a $\delta$-random resolution refutation
distribution of clauses
$Clique_{n,\omega} \cup Color_{n, \xi}$ having $k$ steps.
For a sample $\bs$ define the set $Bad_\bs \subseteq U \times V$
to be the set of all pairs $(u,v) \in U \times V$ such that the
assignment $w(u,v)$ falsifies some clause in $\Delta_\bs$.
An averaging argument implies the following statement.

\begin{lem} 
There exists sample $\bs$ such that $|Bad_\bs| < \delta |U \times V|$.  
\end{lem}

Fix for the rest of the paper one such $\bs$.
Denote by $\pro$ the protocol for $F_{\Delta_\bs}$
constructed from
$\pi_\bs$ as described in Theorem \ref{3.2}.
Put $d := |\Delta_\bs|$.

\begin{lem} \label{15.4.a}
There exists $U'\subseteq U$ and $V'\subseteq V$ such that:
\begin{enumerate}
\item $(U' \times V') \cap Bad_\bs = \emptyset$.

\item $|U'| \geq (1- d \delta^{1/2}) |U|$ and 
$|V'| \geq (1- d \delta^{1/2}) |V|$.
\end{enumerate}
\end{lem}

\proof

\medskip
\noindent
{\bf Claim 1:} {\em The set $Bad_\bs$ is a union of
at most $d'$ rectangles $U_i \times V_i \subseteq U \times V$, for $1 \le d'\le d$.}

\smallskip

For a clause $D$ let $False(D)$ is the set of all $(u,v) \in U \times V$ such that $w(u,v)$ falsifies
$D$. We have that 
$$
Bad_\bs = \bigcup_{D \in \Delta_\bs} False(D)\ .
$$
But for each of at most $d$ possible $D$ 
the set $False(D)$ is a combinatorial 
rectangle as
it consists of all pairs $(u,v) \in U \times V$ satisfying two separate conditions for $u$ and $v$:
that $q^u$ makes all $\bq$-literals in $D$ false and that $v, r^v$ makes all
$\bp$- and $\rr$-literals in $D$ false.

\medskip

Let $\mu_i$ be the measure of $U_i \times V_i$ in $U \times V$ (and so $\mu_i < \delta$). 
The following is obvious.

\medskip
\noindent
{\bf Claim 2:} {\em For each $i \le d'$, either $|U_i| \le \mu_i^{1/2}|U|$ or
$|V_i| \le \mu_i^{1/2}|V|$.}

\medskip

We are now ready to prove the lemma.
Consider the following process. For $i = 1, \dots, d'$ delete from $U$ all elements in $U_i$, if
$|U_i| \le \mu_i^{1/2}|U|$, otherwise delete from $V$ all elements of $V_i$. Let $U'$ and $V'$
be what remains of $U$ and $V$, respectively. 
Because we deleted one side of every rectangle
$U_i \times V_i$, all of them have the empty intersection with $U'\times V '$. 

The measure of $U \setminus U'$ in $U$, as well as 
the measure of $V \setminus V '$ in $V$, 
is bounded above by $\sum_{i \le d'} \mu_i^{1/2} < d \delta^{1/2}$.\qed

\begin{lem} \label{15.4.b}
There exists a monotone protocol for $KW^m$ on $U'\times V'$ of size at most
$k + {n \choose 2}$ and of communication complexity at most $O(\log n)$.
\end{lem}
\proof

Take the protocol $\pro$ for $F_{\Delta_\bs}$ described before Lemma \ref{15.4.a}. 
By the definition of the sets $U'$ and $V'$ the multifunction
$F_{\Delta_\bs}$ restricted to $U' \times V'$ is just $KW^m$
(the condition permitting the extra value $\bot$ is never satisfied).
\qed

\bigskip

\noindent
{\bf Proof of Theorem \ref{main}}:

\medskip

The proof of the $n^{\Omega(\xi^{1/2})}$ lower bound from \cite{AloBop}
for monotone circuits separating $U$ from $V$ 
culminates by comparing two quantities
with the sizes of $U$ and $V$, respectively 
(see the elementary presentation in \cite[Sec.4.3]{BopSip}).
The same argument applies also to separations of any $U'\subseteq U$ from
any $V'\subseteq V$ and the resulting lower bound just gets multiplied by
the smaller of the two measures $|U'|/|U|$ and $|V'|/|V|$.

By Lemmas \ref{15.4.a} and \ref{15.4.b} we have two sets $U', V'$ of relative measures at least
$(1 - d\delta^{1/2})$ and a monotone protocol for $KW^m$ on them of the size at most 
$k+{n \choose 2}$ and communication complexity $O(\log n)$. By Theorem \ref{3.2} 
this yields
a monotone circuit separating $U'$ from $V'$ of size $k n^{O(1)}$. Hence it must hold:
$$
k n^{O(1)}\ \geq (1 - d\delta^{1/2}) n^{\Omega(\xi^{1/2})}
$$
which entails the first inequality in Theorem \ref{main}. 
The second follows from the first one
by estimating $d \le k$: if $k \le 1/(2\delta^{1/2})$ then the factor 
$(1 - d\delta^{1/2})$ is at least $1/2$ and the lower bound
$n^{\Omega(\xi^{1/2})}$ follows.
\qed

\bigskip
\noindent
{\bf Acknowledgements:} I thank E.~Je\v r\' abek, P.~Pudl\' ak and N.~Thapen for comments on early versions
of this paper, and the anonymous referee for suggestions improving the presentation.

%
%
%
%
%
%

\end{document}